\documentclass[10pt]{amsart}

\usepackage{enumerate, amsmath, amsfonts, amssymb, amsthm, mathtools, thmtools, wasysym, graphics, graphicx, xcolor, frcursive,xparse,comment,bbm}
\usepackage[all]{xy}
\usepackage[colorlinks=true, pdfstartview=FitV, linkcolor=blue, citecolor=blue, urlcolor=blue]{hyperref}

\usepackage{url, hypcap}
\hypersetup{colorlinks=true, citecolor=darkblue, linkcolor=darkblue}
\usepackage{mathrsfs}

\usepackage{tikz}
\usetikzlibrary{arrows,automata}
\usetikzlibrary{patterns,snakes}
\usetikzlibrary{shapes.geometric,calc}
\tikzstyle{arrow}=[thick, ->, >=stealth]
\tikzset{
	edge/.style={->,> = latex'}
}
\usetikzlibrary{calc,through,backgrounds,shapes,matrix}
\usepackage{graphicx}

\newcommand{\arxiv}[1]{\href{http://arxiv.org/abs/#1}{\texttt{arXiv:#1}}}

\usepackage[letterpaper,margin=1.25in]{geometry}

\definecolor{darkblue}{rgb}{0.0,0,0.7} 
\definecolor{darkgreen}{rgb}{0, .6, 0} 
\definecolor{lightblue}{rgb}{0,135,147}

\definecolor{red}{rgb}{1,0,0}

\definecolor{darkred}{rgb}{0.7,0,0} 
\definecolor{lightgrey}{rgb}{0.7,0.7,0.7} 

\newtheorem{theorem}{Theorem}[section]
\newtheorem{proposition}[theorem]{Proposition}

\theoremstyle{definition}
\newtheorem{definition}[theorem]{Definition}
\newtheorem{example}[theorem]{Example}

\newtheorem{remark}[theorem]{Remark}

\numberwithin{equation}{section}

\definecolor{darkred}{rgb}{0.7,0,0} 
\newcommand{\defn}[1]{{\color{darkred}\emph{#1}}} 

\usepackage[colorinlistoftodos]{todonotes}

\title[Markov chains]{Markov chains through semigroup graph expansions \newline (a survey)}
  
\author[J.~Rhodes]{John Rhodes}
\address[J. Rhodes]{Department of Mathematics, University of California, Berkeley, CA 94720, U.S.A.}
\email{rhodes@math.berkeley.edu, blvdbastille@gmail.com}

\author[A.~Schilling]{Anne Schilling}
\address[A. Schilling]{Department of Mathematics, UC Davis, One Shields Ave., Davis, CA 95616-8633, U.S.A.}
\email{anne@math.ucdavis.edu}

\date{\today}
\keywords{Markov chains, Karnofsky--Rhodes expansion, McCammond expansion, stationary distribution, mixing time}
\subjclass[2010]{Primary 20M30, 60J10; Secondary 20M05, 60B15, 60C05}

\dedicatory{Dedicated to K. S. S. Nambooripad}

\begin{document}

\begin{abstract}
We review the recent approach to Markov chains using the Karnofksy--Rhodes and McCammond
expansions in semigroup theory by the authors and illustrate them by two examples.
\end{abstract}

\maketitle

\section{Introduction}

A Markov chain is a model that describes transitions between states in a state space according to certain probabilistic rules. 
The defining characteristic of a Markov chain is that the transition from one state to another only depends on the current
state and the elapsed time, but not how it arrived there. In other words, a Markov chain is  ``memoryless''. Markov chains
have many applications, from stock performances, population dynamics to traffic models.

In this paper, we consider Markov chains with a finite state space $\Omega$. The Markov chain can be described
pictorially by its \defn{transition diagram}, where the vertices are the states in $\Omega$ and the labelled directed arrows 
between the vertices indicate the transitions. For example
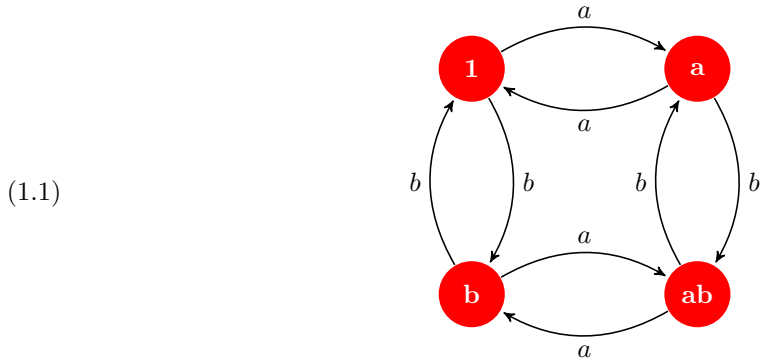
\begin{equation}
\label{equation.transition diagram}
\raisebox{-2.2cm}{
\begin{tikzpicture}[->,>=stealth',shorten >=1pt,auto,node distance=3cm,
                    semithick]
  \tikzstyle{every state}=[fill=red,draw=none,text=white]

  \node[state]         (I) {$\mathbf{1}$};
  \node[state]         (A) [right of=I] {$\mathbf{a}$};
  \node[state]         (B) [below of=I] {$\mathbf{b}$};
  \node[state]         (C) [right of =B] {$\mathbf{ab}$};

  \path (I) edge [bend left] node {$a$} (A)
                 edge [bend left] node {$b$} (B)
           (A) edge[bend left] node {$a$} (I)
                 edge[bend left] node {$b$} (C)
           (B) edge[bend left] node {$b$} (I)
                 edge[bend left] node {$a$} (C)
           (C) edge[bend left] node {$a$} (B)
                 edge[bend left] node {$b$} (A);
\end{tikzpicture}}
\end{equation}
is the transition diagram for a Markov chain with state space $\Omega = \{\mathbf{1},\mathbf{a},\mathbf{b},\mathbf{ab}\}$.
 Associating the probability $0\leqslant x_a \leqslant 1$ to the transition arrows labeled $a$ and $0\leqslant x_b \leqslant 1$ 
to the transition arrows labeled $b$ yields a Markov chain assuming that $x_a+x_b=1$. As the picture indicates, if the Markov 
chain is in state $\mathbf{1}$, then with probability $x_a$ it transitions to state $\mathbf{a}$ and with probability $x_b$ it 
transitions to state $\mathbf{b}$, and so on. Note that every vertex has one outgoing arrow labeled $a$ and one outgoing 
arrow labeled $b$.

The \defn{transition matrix} $T$ of a Markov chain is a matrix of dimension $|\Omega|\times |\Omega|$.
Let $A$ be the set of all edge labels in the transition diagram. Then the entry $T_{s',s}$ in row $s'\in \Omega$ and 
column $s\in \Omega$ in $T$ is
\[
	T_{s',s} = \sum_{\substack{a\in A\\ s \stackrel{a}{\longrightarrow} s'}} x_a,
\]
where the sum is over all $a\in A$ such that $s \stackrel{a}{\longrightarrow} s'$ is an edge in the transition diagram.
For the transition diagram in~\eqref{equation.transition diagram}, the transition matrix is given by
(ordering the states as $(\mathbf{1},\mathbf{a},\mathbf{b},\mathbf{ab})$)
\[
	T = \begin{pmatrix}
	0&x_a&x_b&0\\ x_a&0&0&x_b\\ x_b&0&0&x_a\\ 0&x_b&x_a&0
	\end{pmatrix}.
\]
Note that, since every vertex has precisely one outgoing edge labeled $a\in A$ and since $\sum_{a\in A} x_a=1$, the
column sums of $T$ are equal to one, namely $\sum_{s'\in \Omega} T_{s',s}=1$.
Starting with a distribution $\nu$ (where state $s\in \Omega$ occurs with probability $\nu_s$ and $\sum_{s\in \Omega}
\nu_s=1$), the distribution of states after $t$ steps in the Markov chain is $T^t \nu$.

Two fundamental questions are to find the \defn{stationary distribution} and the \defn{mixing time} of the Markov chain.
Intuitively speaking, the stationary distribution is the distribution of states that the Markov chain will tend to when the chain 
runs for a long time. For ergodic Markov chains (to be defined later), the stationary distribution is unique and is the right
eigenvector of $T$ of eigenvalue one. The mixing time measures how quickly the distribution approaches the stationary 
distribution.

In this paper, we review the approach of~\cite{RhodesSchilling.2019a,RhodesSchilling.2019b} to compute the 
stationary distribution of a finite Markov chain using expansions of semigroups. For example, the Markov chain with
transition diagram~\eqref{equation.transition diagram} can be described using the dihedral group
\[
	D_2 = \langle a,b \mid a^2=b^2=1, (ab)^2=1\rangle,
\]
generated by two reflections $a,b$. This group is isomorphic to $Z_2 \times Z_2$.
State $\mathbf{1}$ corresponds to the identity, state $\mathbf{a}$ corresponds to $a$, 
state $\mathbf{b}$ to $b$, and state $\mathbf{ab}$ corresponds to $ab=ba$. The transition from state $s$ to $s'$ is given 
by left multiplication by one of the generators in $A=\{a,b\}$. In general, it is always possible to describe a finite state Markov 
chain via a semigroup~\cite{LevinPeres.2017, ASST.2015} by the \defn{random letter representation}.

\begin{figure}[t]
\begin{center}
\begin{tikzpicture}[auto]
\node (A) at (0, 0) {$\mathbbm{1}$};
\node (B) at (-2,-1) {$\mathbf{a}$};
\node(C) at (2,-1) {$\mathbf{b}$};
\node(D) at (0,-2) {$\mathbf{ab}$};
\node(E) at (0,-3) {$\mathbf{1}$};
\draw[edge,blue,thick] (A) -- (B) node[midway, above] {$a$};
\draw[edge,blue,thick] (A) -- (C) node[midway, above] {$b$};
\draw[edge,thick] (B) -- (D) node[midway, above] {$b$};
\draw[edge,thick] (D) -- (B);
\draw[edge,thick] (C) -- (D) node[midway, above] {$a$};
\draw[edge,thick] (D) -- (C);
\draw[edge,thick] (B) -- (E) node[midway, below] {$a$};
\draw[edge,thick] (E) -- (B);
\draw[edge,thick] (C) -- (E) node[midway, below] {$b$};
\draw[edge,thick] (E) -- (C);
\end{tikzpicture}
\end{center}
\caption{\label{figure.right cayley}The right Cayley graph $\mathsf{RCay}(D_2,\{a,b\})$ of the dihedral group
 with generators $A=\{a,b\}$. Transition edges are indicated in blue. Double edges mean that 
 right multiplication by the label for either vertex yields the other vertex.}
\end{figure}
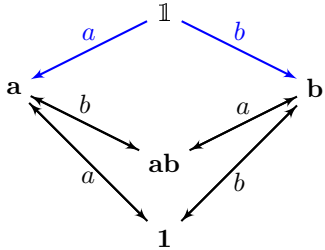

Our focus in this paper is to compute the stationary distribution from the McCammond and Karnofsky--Rhodes expansion 
of the right Cayley graph of the underlying semigroup $S$ with generators in $A$. The right Cayley graph 
$\mathsf{RCay}(S,A)$ of the dihedral group $S=D_2$ with generators $A=\{a,b\}$ is depicted in 
Figure~\ref{figure.right cayley}. In Section~\ref{section.stationary}, we will define the right Cayley graph of 
a semigroup, introduce its Karnofsky--Rhodes and McCammond expansion, and review the main results
from~\cite{RhodesSchilling.2019a,RhodesSchilling.2019b} on how to compute the stationary distribution of the
Markov chain from these. We will illustrate the results in terms of two (running) examples.

\subsection*{Acknowledgments}
We would like to thank Igor Pak for discussions.

The last author was partially supported by NSF grants DMS--1760329 and DMS--1764153. 
This material is based upon work supported by the Swedish Research Council under
grant no. 2016-06596 while the last author was in residence at Institut 
Mittag--Leffler in Djursholm, Sweden during Spring~2020.
The authors thank the organizers of the International conference on semigroups and applications
held at Cochin University of Science and Technology in India December 9-12, 2019, where this work
was presented.

\section{Stationary distribution from semigroup expansions}
\label{section.stationary}

\subsection{Markov chains}
\label{section.markov}

Let $\mathcal{M}$ be a finite Markov chain with state space $\Omega$ and transition matrix $T$.
A Markov chain is \defn{irreducible} if the transition diagram of the Markov chain is strongly connected. It is \defn{aperiodic}
if the greatest common divisor of the cycle lengths in the transition diagram of the Markov chain is one. Furthermore, a 
Markov chain is \defn{ergodic} if it is both irreducible and aperiodic. By the Perron--Frobenius Theorem, an ergodic Markov 
chain has a unique stationary distribution $\Psi$ and $T^t \nu$ converges to $\Psi$ as $t\to \infty$ for any initial state $\nu$. 
In fact, the \defn{stationary distribution} is the right eigenvector of eigenvalue one of $T$
\[
	T \Psi = \Psi.
\]

An important question is how quickly does the Markov chain converge to the stationary distribution.
In Markov chain theory, distance is usually the total variation distance. The total variation distance between two probability 
distributions $\nu$ and $\mu$ is defined as
\[
	\|\nu - \mu\| = \max_{A \subseteq \Omega} |\nu(A) - \mu(A)|.
\]
For a given small $\epsilon>0$, the \defn{mixing time} $t_\mathsf{mix}$ is the smallest $t$ such that
\[
	\| T^t \nu - \Psi \| \leqslant \epsilon.
\]

Brown and Diaconis~\cite{BrownDiaconis.1998,Brown.2000} analyzed Markov chains associated to left regular bands.
In particular, they showed~\cite[Theorem 0]{Brown.2000} that the total variational distance from stationarity after $t$ steps 
is bounded above by the probability $\mathsf{Pr}(\tau> t)$, where $\tau$ is the first time that the walk hits a certain ideal. 
The arguments in Brown and Diaconis~\cite{BrownDiaconis.1998} were generalized to Markov chains for
$\mathscr{R}$-trivial semigroups~\cite{ASST.2015}. A unified theory for Markov chains for any finite semigroup
was developed in~\cite{RhodesSchilling.2019a,RhodesSchilling.2019b}.

As explained in~\cite[Proposition 1.5]{LevinPeres.2017} and~\cite[Theorem 2.3]{ASST.2015}, every finite state Markov chain 
$\mathcal{M}$ has a random letter representation, that is, a representation of a semigroup $S$ acting on the left on the 
state space $\Omega$. In this setting, we transition $s \stackrel{a}{\longrightarrow} s'$ with probability 
$0\leqslant x_a\leqslant 1$, where $s, s'\in \Omega$, $a\in S$ and $s'=a.s$ is the action of $a$ on the state $s$. 
It is enough to consider the semigroup $S$ generated by the elements $a\in A$ with $x_a>0$.

A two-sided \defn{ideal} $I$ (or ideal for short) is a subset $I \subseteq S$ such that $u I v \subseteq I$ for all 
$u,v \in S^{\mathbbm{1}}$, where $S^{\mathbbm{1}}$ is the semigroup $S$ with the identity $\mathbbm{1}$ added (even if
$S$ already contains an identity). If $I,J$ are ideals of $S$, then $IJ \subseteq I \cap J$, so that $I \cap J \neq \emptyset$. 
Hence every finite semigroup has a unique \defn{minimal ideal} denoted $K(S)$.
An ideal $K(S)$ is \defn{left zero} if $xy=x$ for all $x,y\in K(S)$.

We will determine the stationary distribution from certain expansions of the right Cayley graph $\mathsf{RCay}(S,A)$
of the underlying semigroup $S$ with generators $A$. The Markov chain itself is a random walk on the minimal ideal
$K(S)$ by the left action.

\subsection{Right Cayley graphs}
We begin with the definition of a graph.

\begin{definition}[Graph]
A \defn{labeled directed graph} $\Gamma$ (or \defn{graph} for short) consists of a vertex set $V(\Gamma)$, an edge 
set $E(\Gamma)$, and a labelling set $A$. An edge $e\in E(\Gamma)$ is a tuple $e=(v,a,w) \in V(\Gamma)\times A \times
V(\Gamma)$. We often also write $e\colon v \stackrel{a}{\longrightarrow} w$.
\end{definition} 

A \defn{path} $p$ from vertex $v$ to vertex $w$ in a graph $\Gamma$ is a sequence of edges
\[
	p = \left(v = v_0 \stackrel{a_1}{\longrightarrow} v_1 \stackrel{a_2}{\longrightarrow} \cdots 
	\stackrel{a_\ell}{\longrightarrow} v_\ell =w\right),
\]
where each tuple $(v_i,a_{i+1},v_{i+1}) \in E(\Gamma)$ for $0\leqslant i<\ell$. The initial (resp. terminal) vertex $v$ 
(resp. $w$) of $p$ is denoted by $\iota(p)$ (resp. $\tau(p)$). The \defn{length} of $p$ is $\ell(p):= \ell$
and $a_1 \ldots a_\ell$ is called the label of the path. 

We can define a preorder $\prec$ on $V(\Gamma)$ by $v \prec w$ if there is a path from $v$ to $w$ in $\Gamma$.
This induces an equivalence relation $\sim$ on $V(\Gamma)$, where $v \sim w$ if $v \prec w$ and $w \prec v$.
A \defn{strongly connected component} of $\Gamma$ is a $\sim$-equivalence class.

\begin{definition}[Rooted graph]
A \defn{rooted graph} is a pair $(\Gamma,r)$, where $\Gamma$ is a graph and $r \in V(\Gamma)$, such that
$r \prec v$ for all $v \in V(\Gamma)$.
\end{definition}

A path is called \defn{simple} if it visits no vertex twice. Empty (or trivial) paths are considered simple.
For a rooted graph $(\Gamma,r)$, let $\mathsf{Simple}(\Gamma,r)$ be the set of simple paths of $\Gamma$
starting at $r$ (including the empty path).

\begin{definition}[Right Cayley graph]
Let $(S,A)$ be a finite semigroup $S$ together with a set of generators $A$. 
The \defn{right Cayley graph} $\mathsf{RCay}(S,A)$ of $S$ with respect to $A$ is the rooted graph with
vertex set $V(\mathsf{RCay}(S,A)) = S^{\mathbbm{1}}$, root $r=\mathbbm{1} \in S^{\mathbbm{1}}$, and edges 
$s \stackrel{a}{\longrightarrow} s'$ for all $(s,a,s') \in S^{\mathbbm{1}} \times A \times S^{\mathbbm{1}}$, such
that $s'=sa$ in $S^{\mathbbm{1}}$.  
\end{definition}

An example of a right Cayley graph is given in Figure~\ref{figure.right cayley}.

For a semigroup $S$, two elements $s,s'\in S$ are in the same $\mathscr{R}$-class if the corresponding right ideals 
are equal, that is, $s S^{\mathbbm{1}} = s'S^{\mathbbm{1}}$. The strongly connected components of $\mathsf{RCay}(S,A)$ 
are precisely the $\mathscr{R}$-classes of $S^{\mathbbm{1}}$. In other words, the vertices of a strongly connected 
component are exactly the vertices that represent the elements in an $\mathscr{R}$-class of $S^{\mathbbm{1}}$. 
Edges that go between distinct strongly connected components will turn out to play an important role in the 
Karnofksy--Rhodes expansion.

\begin{definition}[Transition edges]
Let $\Gamma$ be a graph. Then $e=(v,a,w) \in E(\Gamma)$ with $v,w\in V(\Gamma)$ and $a\in A$ is 
a \defn{transition edge} if $v \not \sim w$. In other words, there is no path from $w$ to $v$ in $\Gamma$.
\end{definition}

In Figure~\ref{figure.right cayley}, the transition edges are indicated in blue. Note that the edges leaving $\mathbbm{1}$
in the right Cayley graph are always transitional. Other edges might or might not be transitional. In this example
$K(S)$ consists of all vertices in $\mathsf{RCay}(S,A)$ except the root $\mathbbm{1}$.

\subsection{The Karnofsky--Rhodes expansion}
\label{section.KR}

To compute explicit expressions for the stationary distributions of Markov chains on finite semigroups, we need the
\defn{Karnofsky--Rhodes expansion}~\cite{Elston.1999} of the right Cayley graph $\mathsf{RCay}(S,A)$.
See also~\cite[Definition~4.15]{MRS.2011} and~\cite[Section 3.4]{MSS.2015}. 

Denote by $(A^+,A)$ the free semigroup with generators in $A$. In other words, $A^+$ is the set of all words 
$a_1 \ldots a_\ell$ of length $\ell \geqslant 1$ over $A$ with multiplication given by concatenation. 
Furthermore, let $A^\star = A^+ \cup \{1\}$, so that $A^\star$ is $A^+$ with the identity added; it is the free monoid 
generated by $A$. 

\begin{definition}[Karnofksy--Rhodes expansion]
The \defn{Karnofsky--Rhodes expansion} $\mathsf{KR}(S,A)$ is obtained as follows. Start with the right Cayley graph 
$\mathsf{RCay}(A^+,A)$. Identify the endpoints of two paths in $\mathsf{RCay}(A^+,A)$ 
\begin{equation*}
	p := \left( \mathbbm{1} \stackrel{a_1}{\longrightarrow} v_1 \stackrel{a_2}{\longrightarrow} \cdots \stackrel{a_\ell}
	{\longrightarrow} v_\ell \right)
	\quad \text{and} \quad 
	p' := \left( \mathbbm{1} \stackrel{a'_1}{\longrightarrow} v'_1 \stackrel{a'_2}{\longrightarrow} \cdots 
	\stackrel{a'_{\ell'}}{\longrightarrow} v'_{\ell'} \right)
\end{equation*}
in $\mathsf{KR}(S,A)$ if and only if the corresponding paths in $\mathsf{RCay}(S,A)$
\begin{equation*}
	[p]_S := \left( \mathbbm{1} \stackrel{a_1}{\longrightarrow} [w_1]_S \stackrel{a_2}{\longrightarrow} \cdots 
	\stackrel{a_\ell}{\longrightarrow} [w_\ell]_S \right)
	\quad \text{and} \quad 
	[p']_S := \left( \mathbbm{1} \stackrel{a'_1}{\longrightarrow} [w'_1]_S \stackrel{a'_2}{\longrightarrow} \cdots 
	\stackrel{a'_{\ell'}}{\longrightarrow} [w'_{\ell'}]_S \right),
\end{equation*}
where $w_i=a_1 a_2 \ldots a_i$ and $w_i' = a_1' a_2' \ldots a'_i$, end at the same vertex $[w_\ell]_S = [w'_{\ell'}]_S$ 
and in addition the set of transition edges of $[p]_S$ and $[p']_S$ in $\mathsf{RCay}(S,A)$ is equal. 
\end{definition}

\begin{figure}[t]
\begin{center}
\begin{tikzpicture}[auto]
\node (A) at (0, 0) {$\mathbbm{1}$};
\node (B) at (-1,-1) {$a$};
\node(C) at (1,-1) {$b$};
\node(D) at (-1,-2) {$ab$};
\node(E) at (1,-2) {$ba$};
\node(F) at (-2.5,-2) {$a^2$};
\node(G) at (2.5,-2) {$b^2$};
\node(H) at (-1,-3) {$a^2b=aba$};
\node(I) at (1,-3) {$bab=b^2a$};
\draw[edge,blue,thick] (A) -- (B) node[midway, above] {$a$\;};
\draw[edge,thick,blue] (A) -- (C) node[midway, above] {\;$b$};
\draw[edge,thick] (B) -- (F) node[midway,left] {$a$\;};
\draw[edge,thick] (F) -- (B);
\draw[edge,thick] (B) -- (D) node[midway,right] {$b$};
\draw[edge,thick] (D) -- (B);
\draw[edge,thick] (F) -- (H) node[midway,left] {$b$\;};
\draw[edge,thick] (H) -- (F);
\draw[edge,thick] (D) -- (H) node[midway,right] {\;$a$};
\draw[edge,thick] (H) -- (D);
\draw[edge,thick] (C) -- (E) node[midway,left] {$a$\;};
\draw[edge,thick] (E) -- (C);
\draw[edge,thick] (E) -- (I) node[midway,left] {$b$\;};
\draw[edge,thick] (I) -- (E);
\draw[edge,thick] (C) -- (G) node[midway,right] {\;$b$};
\draw[edge,thick] (G) -- (C);
\draw[edge,thick] (G) -- (I) node[midway,right] {\;$a$};
\draw[edge,thick] (I) -- (G);
\end{tikzpicture}
\end{center}
\caption{\label{figure.KR} The Karnofsky--Rhodes expansion $\mathsf{KR}(S,A)$ of the right Cayley graph of
Figure~\ref{figure.right cayley}.}
\end{figure}
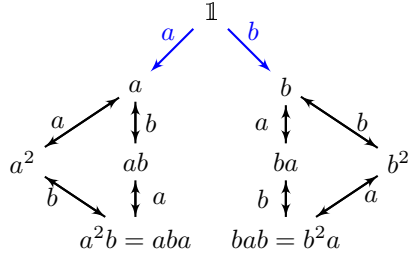

An example for $\mathsf{KR}(S,A)$ is given in Figure~\ref{figure.KR}. In this figure, the paths $a^2b$ and $aba$
are equal because they end in the same vertex when projected onto $S$ and they share the same transition edge, 
which is the first $a$. On the other hand, the paths $ab$ and $ba$ are distinct even though $ab=ba$ in $D_2$
because for the first path the transition edge is the first $a$ and for the second path the transition edge is the first $b$.

\begin{proposition}\cite[Proposition 2.15]{RhodesSchilling.2019a}
\label{proposition.KR Cayley}
$\mathsf{KR}(S,A)$ is the right Cayley graph of a semigroup, also denoted by $\mathsf{KR}(S,A)$.
\end{proposition}

\subsection{The McCammond expansion}
\label{section.mccammond}

The McCammond expansion~\cite{MRS.2011} of a rooted graph is intimately related to the 
unique simple path property.

\begin{definition}[Unique simple path property]
\label{definition.unique simple path}
A rooted graph $(\Gamma,r)$ has the \defn{unique simple path property} if for each vertex $v\in V(\Gamma)$ there is
a unique simple path from the root $r$ to $v$.
\end{definition}

As proven in~\cite[Proposition 2.32]{MRS.2011}, the unique simple path property is equivalent to $(\Gamma,r)$
admitting a unique directed spanning tree $\mathsf{T}$. Note that the unique simple path property not only depends 
on the graph $\Gamma$, but also on the chosen root $r$. In this paper, we always choose $r= \mathbbm{1}$.
It was established in~\cite[Section 2.7]{MRS.2011} that every rooted graph $(\Gamma,r)$ has a universal simple
cover, which has the unique simple path property.

If $p$ and $q$ are paths, $\ell(q)=k \leqslant \ell(p)$, and the first 
$k+1$ vertices and $k$ edges of $p$ and $q$ agree, we say that $q$ is an \defn{initial segment} of $p$, written 
$q \subseteq p$.

\begin{definition}[McCammond expansion]
\label{definition.mccammond}
For a rooted graph $(\Gamma,r)$, define its \defn{McCammond expansion} $(\Gamma^{\mathsf{Mc}},r)$ as the graph with
\begin{equation*}
\begin{split}
	V(\Gamma^{\mathsf{Mc}}) &= \mathsf{Simple}(\Gamma,r),\\
	E(\Gamma^{\mathsf{Mc}}) &= \{(p,a,q) \in V(\Gamma^{\mathsf{Mc}}) \times A \times V(\Gamma^{\mathsf{Mc}}) \mid
	(\tau(p), a, \tau(q)) \in E(\Gamma),\\
	& \qquad \qquad  \qquad \ell(q) = \ell(p)+1 \text{ or } (q\subseteq p \text{ and } \ell(q)\leqslant \ell(p))\}.
\end{split}
\end{equation*}
\end{definition}

Note that by definition there are two types of edges $(p,a,q) \in E(\Gamma^{\mathsf{Mc}})$: either $\ell(q)=\ell(p)+1$
or $\ell(q) \leqslant \ell(p)$ as paths in $\mathsf{Simple}(\Gamma,r)$. The spanning tree $\mathsf{T}$ has vertex set
$V(\Gamma^{\mathsf{Mc}})$ and only those edges $(p,a,q) \in E(\Gamma^{\mathsf{Mc}})$ such that $\ell(q)=\ell(p)+1$.

\begin{figure}[t]
\begin{center}
\begin{tikzpicture}[auto]
\node (A) at (0, 0) {$\mathbbm{1}$};
\node (B) at (-3,-1) {$a$};
\node(C) at (3,-1) {$b$};
\node(D) at (-2,-2) {$ab$};
\node(E) at (2,-2) {$ba$};
\node(F) at (-4,-2) {$a^2$};
\node(FF) at (-4,-3.5) {$a^2b$};
\node(FFF) at (-4,-5) {$a^2ba$};
\node(G) at (4,-2) {$b^2$};
\node(H) at (-2,-3.5) {$aba$};
\node(HH) at (-2,-5) {$abab$};
\node(EE) at (2,-3.5) {$bab$};
\node(EEE) at (2,-5) {$baba$};
\node(GG) at (4,-3.5) {$b^2a$};
\node(GGG) at (4,-5) {$b^2ab$};
\draw[edge,blue,thick] (A) -- (B) node[midway, above] {$a$\;};
\draw[edge,thick,blue] (A) -- (C) node[midway, above] {\;$b$};
\draw[edge,thick] (B) -- (F) node[midway,right] {$a$};
\path (F) edge[->,thick, red,dashed, bend left=30] node[midway,left] {$a$} (B);
\draw[edge,thick] (F) -- (FF) node[midway,right] {$b$};
\path (FF) edge[->,thick, red, dashed, bend left=30] node[midway,left] {$b$} (F);
\draw[edge,thick] (FF) -- (FFF) node[midway,right] {$a$};
\path (FFF) edge[->,thick, red, dashed, bend left=30] node[midway,left] {$a$} (FF);
\path (FFF) edge[->,thick, red,dashed, bend left=90] node[midway,left] {$b$} (B);
\draw[edge,thick] (B) -- (D) node[midway,left] {$b$};
\path (D) edge[->,thick, red, dashed, bend right=30] node[midway,right] {$b$} (B);
\draw[edge,thick] (D) -- (H) node[midway,left] {$a$};
\path (H) edge[->,thick, red, dashed, bend right=30] node[midway,right] {$a$} (D);
\draw[edge,thick] (H) -- (HH) node[midway,left] {$b$};
\path (HH) edge[->,thick, red, dashed, bend right=30] node[midway,right] {$b$} (H);
\path (HH) edge[->,thick, red,dashed, bend right=90] node[midway,right] {$a$} (B);
\draw[edge,thick] (C) -- (E) node[midway,right] {$a$};
\path (E) edge[->,thick, red,dashed, bend left=30] node[midway,left] {$a$} (C);
\draw[edge,thick] (E) -- (EE) node[midway,right] {$b$};
\path (EE) edge[->,thick, red,dashed, bend left=30] node[midway,left] {$b$} (E);
\draw[edge,thick] (EE) -- (EEE) node[midway,right] {$a$};
\path (EEE) edge[->,thick, red,dashed, bend left=30] node[midway,left] {$a$} (EE);
\draw[edge,thick] (C) -- (G) node[midway,left] {$b$};
\path (G) edge[->,thick, red,dashed, bend right=30] node[midway,right] {$b$} (C);
\draw[edge,thick] (G) -- (GG) node[midway,left] {$a$};
\path (GG) edge[->,thick, red,dashed, bend right=30] node[midway,right] {$a$} (G);
\draw[edge,thick] (GG) -- (GGG) node[midway,left] {$b$};
\path (GGG) edge[->,thick, red,dashed, bend right=30] node[midway,right] {$b$} (GG);
\path (EEE) edge[->,thick, red,dashed, bend left=90] node[midway,left] {$b$} (C);
\path (GGG) edge[->,thick, red,dashed, bend right=90] node[midway,right] {$a$} (C);
\end{tikzpicture}
\end{center}
\caption{\label{figure.mccammond} The McCammond expansion $\mathsf{Mc} \circ \mathsf{KR}(D_2,\{a,b\})$ 
of Figure~\ref{figure.KR}.}
\end{figure}
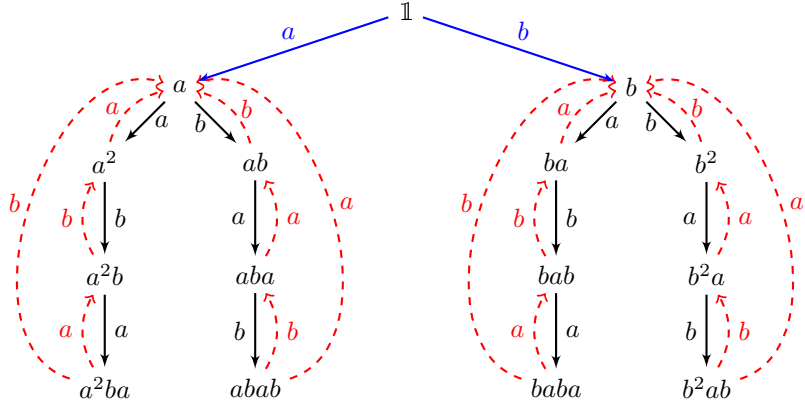

From now on choose $r = \mathbbm{1}$. The simple path 
\[
	\mathbbm{1} \stackrel{a_1}{\longrightarrow} v_1 \stackrel{a_2}{\longrightarrow} \cdots \stackrel{a_\ell}{\longrightarrow} 
	v_\ell
\]
in $\mathsf{Simple}(\Gamma,\mathbbm{1})$ is naturally indexed by the word $a_1 a_2 \ldots a_\ell$. We will use this labeling
for the McCammond expansion of $\mathsf{KR}(S,A)$. In particular, if $a_1 a_2 \ldots a_\ell \in 
\mathsf{Simple}(\Gamma,\mathbbm{1})$ and $a_1 a_2 \ldots a_\ell a \in \mathsf{Simple}(\Gamma,\mathbbm{1})$, then 
the edge $a_1 a_2 \ldots a_\ell \stackrel{a}{\longrightarrow} a_1 a_2 \ldots a_\ell a$ is in the spanning tree $\mathsf{T}$.
Otherwise we have $a_1 a_2 \ldots a_\ell \stackrel{a}{\longrightarrow} a_1 a_2 \ldots a_k$ for some unique $1\leqslant k < \ell$.
Thus under the right action of $a \in A$ on $a_1 a_2 \ldots a_\ell$, we either move forward in the spanning tree or
fall backwards somewhere on the unique geodesic from $\mathbbm{1}$ to $a_1 a_2 \ldots a_\ell$, 
but staying in the same $\mathscr{R}$-class. An example of a McCammond expansion of a Karnofsky--Rhodes graph is 
given in Figure~\ref{figure.mccammond}.

For a non-simple path in $(\Gamma^{\mathsf{Mc}},\mathbbm{1})$, we can remove loops; it does not matter in which order
these loops are removed. This is also known as the Church--Rosser property~\cite{Church.Rosser.1936} or a
Knuth--Bendix rewriting system. This is proved in~\cite{MRS.2011}.

We denote the McCammond expansion of a semigroup $(S,A)$ with generators in $A$ by $\mathsf{Mc}(S,A)$, which is 
the McCammond expansion of its right Cayley graph. 

\subsection{Stationary distribution}

Denote by $\mathcal{M}(S,A)$ the Markov chain associated with the semigroup $(S,A)$. As mentioned before,
this is the random walk on $K(S)$ by the left action. The probability $x_a$ is associated with generator $a\in A$.
The stationary distribution for $\mathcal{M}(S,A)$ was computed in~\cite{RhodesSchilling.2019a} using
$\mathsf{Mc}\circ \mathsf{KR}(S,A)$. The treatment depends on whether the minimal ideal
$K(S)$ is left zero or not. We first start with the former case, stated in~\cite[Corollaries 2.23 \& 2.28]{RhodesSchilling.2019a}.

\begin{theorem}
\label{theorem.stationary}
If $K(S)$ is left zero, the stationary distribution of the Markov chain $\mathcal{M}(S,A)$ labeled by $w \in K(S)$
is given by
\begin{equation*}
	\Psi^{\mathcal{M}(S,A)}_w 
	= \sum_{p} \; \prod_{a\in p} x_a,
\end{equation*}
where the sum is over all paths $p$ in $\mathsf{Mc}\circ\mathsf{KR}(S,A)$ starting at $\mathbbm{1}$ and ending
in $s$ such that $[s]_S=w$.
\end{theorem}

The case when $K(S)$ is not left zero was treated in~\cite[Section 2.9]{RhodesSchilling.2019a} by adding a zero
element $\square$ to the semigroup $S$ and the generators $A$. This new generator $\square$ has its own
probability $x_\square$. The minimal ideal of the semigroup $(S\cup \{\square\}, A \cup \{\square\})$ is left zero
and by taking the limit $x_\square \to 0$, the stationary distribution of the original Markov chain 
$\mathcal(S,A)$ is retrieved as stated in~\cite[Corollary~2.33]{RhodesSchilling.2019a}.

\begin{theorem}
\label{theorem.stationary general}
If $K(S)$ is not left zero, the stationary distribution of the Markov chain $\mathcal{M}(S,A)$ labeled by $w \in K(S)$
is given by
\begin{equation}
\label{equation.psi limit}
	\Psi^{\mathcal{M}(S,A)}_w 
	= \lim_{x_\square \to 0} \bigl(\sum_{p}\; \prod_{a\in p} x_a \bigr),
\end{equation}
where the sum is over all paths p in $\mathsf{Mc}\circ\mathsf{KR}(S \cup \{\square\},A \cup \{\square\})$
starting at $\mathbbm{1}$ and ending in $s$ such that $[s]_{S\cup \{\square\}}=w \square$.
\end{theorem}

In~\cite{RhodesSchilling.2019b}, we developed a strategy using \defn{loop graphs} to compute the expressions
in Theorems~\ref{theorem.stationary} and~\ref{theorem.stationary general} as rational functions in the probabilities $x_a$ 
for $a\in A$. This is done in two steps:
\begin{enumerate}
\item Using McCammond's  $\mathsf{Pict}$ map, we map the McCammond expansion together with a simple
path from $\mathbbm{1}$ to the ideal to a loop graph. See Section~\ref{section.loop}.
\item The set of all paths from $\mathbbm{1}$ to the element in the ideal in a loop graph can be written as a 
Kleene expression. The Kleene expression immediately yields a rational expression for the stationary distribution.
See Section~\ref{section.kleene}.
\end{enumerate}

\subsection{Loop graphs}
\label{section.loop}

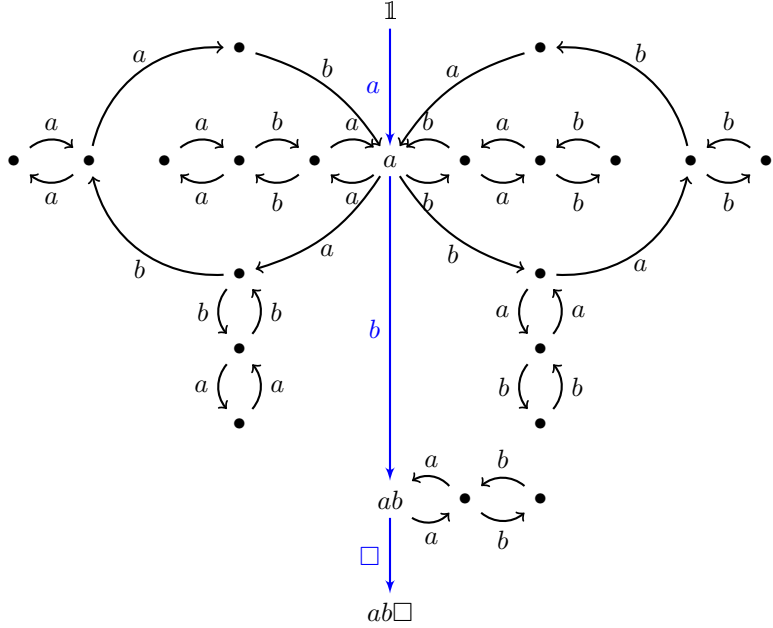
\begin{figure}
\begin{tikzpicture}[auto]
\node (A) at (0, 0) {$\mathbbm{1}$};
\node (B) at (0, -2) {$a$};
\node (C) at (0, -6.5) {$ab$};
\node (D) at (0, -8) {$ab\square$};
\draw[edge,blue,thick] (A) -- (B) node[midway, left] {$a$};
\draw[edge,blue,thick] (B) -- (C) node[midway, left] {$b$};
\draw[edge,blue,thick] (C) -- (D) node[midway, left] {$\square$};
\node (P1) at (1,-6.5) {$\bullet$};
\node (P2) at (2,-6.5) {$\bullet$};
\path (C) edge[->,thick,bend right=40] node[midway,below] {$a$} (P1);
\path (P1) edge[->,thick,bend right=40] node[midway,above] {$a$} (C);
\path (P1) edge[->,thick,bend right=40] node[midway,below] {$b$} (P2);
\path (P2) edge[->,thick,bend right=40] node[midway,above] {$b$} (P1);
\node (L1) at (-1,-2) {$\bullet$};
\node (L2) at (-2,-2) {$\bullet$};
\node (L3) at (-3,-2) {$\bullet$};
\path (B) edge[->,thick,bend left=40] node[midway,below] {$a$} (L1);
\path (L1) edge[->,thick,bend left=40] node[midway,above] {$a$} (B);
\path (L1) edge[->,thick,bend left=40] node[midway,below] {$b$} (L2);
\path (L2) edge[->,thick,bend left=40] node[midway,above] {$b$} (L1);
\path (L2) edge[->,thick,bend left=40] node[midway,below] {$a$} (L3);
\path (L3) edge[->,thick,bend left=40] node[midway,above] {$a$} (L2);
\node (R1) at (1,-2) {$\bullet$};
\node (R2) at (2,-2) {$\bullet$};
\node (R3) at (3,-2) {$\bullet$};
\path (B) edge[->,thick,bend right=40] node[midway,below] {$b$} (R1);
\path (R1) edge[->,thick,bend right=40] node[midway,above] {$b$} (B);
\path (R1) edge[->,thick,bend right=40] node[midway,below] {$a$} (R2);
\path (R2) edge[->,thick,bend right=40] node[midway,above] {$a$} (R1);
\path (R2) edge[->,thick,bend right=40] node[midway,below] {$b$} (R3);
\path (R3) edge[->,thick,bend right=40] node[midway,above] {$b$} (R2);
\node (BL1) at (-2,-3.5) {$\bullet$};
\node (BL2) at (-4,-2) {$\bullet$};
\node (BL3) at (-2,-0.5) {$\bullet$};
\path (B) edge[->,thick,bend left=20] node[midway,below] {$a$} (BL1);
\path (BL1) edge[->,thick,bend left=40] node[midway,below] {$b$} (BL2);
\path (BL2) edge[->,thick,bend left=40] node[midway,above] {$a$} (BL3);
\path (BL3) edge[->,thick,bend left=20] node[midway,above] {$b$} (B);
\node (BR1) at (2,-3.5) {$\bullet$};
\node (BR2) at (4,-2) {$\bullet$};
\node (BR3) at (2,-0.5) {$\bullet$};
\path (B) edge[->,thick,bend right=20] node[midway,below] {$b$} (BR1);
\path (BR1) edge[->,thick,bend right=40] node[midway,below] {$a$} (BR2);
\path (BR2) edge[->,thick,bend right=40] node[midway,above] {$b$} (BR3);
\path (BR3) edge[->,thick,bend right=20] node[midway,above] {$a$} (B);
\node (AL) at (-5,-2) {$\bullet$};
\path (BL2) edge[->,thick,bend left=40] node[midway,below] {$a$} (AL);
\path (AL) edge[->,thick,bend left=40] node[midway,above] {$a$} (BL2);
\node (AR) at (5,-2) {$\bullet$};
\path (BR2) edge[->,thick,bend right=40] node[midway,below] {$b$} (AR);
\path (AR) edge[->,thick,bend right=40] node[midway,above] {$b$} (BR2);
\node (U1) at (-2,-4.5) {$\bullet$};
\path (BL1) edge[->,thick,bend right=40] node[midway,left] {$b$} (U1);
\path (U1) edge[->,thick,bend right=40] node[midway,right] {$b$} (BL1);
\node (U2) at (-2,-5.5) {$\bullet$};
\path (U1) edge[->,thick,bend right=40] node[midway,left] {$a$} (U2);
\path (U2) edge[->,thick,bend right=40] node[midway,right] {$a$} (U1);
\node (U3) at (2,-4.5) {$\bullet$};
\path (BR1) edge[->,thick,bend right=40] node[midway,left] {$a$} (U3);
\path (U3) edge[->,thick,bend right=40] node[midway,right] {$a$} (BR1);
\node (U4) at (2,-5.5) {$\bullet$};
\path (U3) edge[->,thick,bend right=40] node[midway,left] {$b$} (U4);
\path (U4) edge[->,thick,bend right=40] node[midway,right] {$b$} (U3);
\end{tikzpicture}
\caption{Loop graph for $(\mathsf{Mc}\circ \mathsf{KR}(D_2 \cup \{\square\},\{a,b,\square\}), ab\square)$
\label{figure.loop}}
\end{figure}

A \defn{loop} of size $\ell$ is a connected directed graph with $\ell$ vertices such that each vertex has exactly one incoming and 
one outgoing edge. In other words, a loop is a directed circle of $\ell$ vertices. A \defn{loop graph} can be defined recursively.
Start with a directed straight line path. Recursively, attach a loop of an arbitrary finite size to any existing chosen vertex.
Repeat or stop. The edges of the loop graph can be labeled. An example of a loop graph is given in
Figure~\ref{figure.loop}.

Recall that an important property of the McCammond expansion is the unique simple path property
(see Definition~\ref{definition.unique simple path}).
We now define the map $\mathsf{Pict}$ from the set of tuples $(\Gamma,p)$, where $\Gamma$ is a graph
with the unique simple path property and $p$ is a simple path in $\Gamma$ starting at $\mathbbm{1}$, to the set of loop 
graphs. The straight line, that the loop graph is based on, corresponds to the chosen simple path $p$. 
We follow~\cite[Section 3.2]{RhodesSchilling.2019b}.

\begin{definition}[McCammond] \cite[Definition 3.5]{RhodesSchilling.2019b}
\label{definition.pict}
Let $\Gamma$ be a graph with the unique simple path property and $p$ a simple path in $\Gamma$ starting at
$\mathbbm{1}$. Then $\mathsf{Pict}(\Gamma,p)$ is defined by the principle of induction.

\noindent
\textbf{Induction basis:} Set $P=p$ and start at vertex $v_0=\mathbbm{1}$.

\noindent
\textbf{Induction step:} Suppose one is at vertex $v_0 \neq \tau(p)$ on path $p$. Take the edge $e$ from $v_0$ to $v_1$ 
in $p$.
\begin{enumerate}
\item If there is no edge in $\Gamma$ coming into $v_1$ besides $e$, continue with the unique next vertex in $p$,
now denoted $v_1$ (with the current vertex $v_1$ relabeled $v_0$), unless
$v_1=\tau(p)$. If $v_1=\tau(p)$, then output $\mathsf{Pict}(\Gamma,p)=P$.
\item Otherwise there is at least one edge $e' \neq e$ in $\Gamma$ going into $v_1$, given by
$e'=\left(v' \stackrel{a}{\longrightarrow} v_1\right)$ for some $a\in A$. Since $\Gamma$ has the unique simple path property by 
assumption, there must be a unique simple path starting at $\mathbbm{1}$ going to $v_0$ along the path $p$ followed by 
the path $p'$ starting at $v_0$, going along $e$ to $v_1$, and ending at $v'$.
\begin{enumerate}
\item Run the induction on $p'$ in a subgraph $\Gamma'$ of $\Gamma$, 
consisting of all edges and vertices on circuits containing a vertex of $p'$.
Note that $p'$ is simple in $\Gamma'$. The output is $P'=\mathsf{Pict}(\Gamma',p')$.
\item Modify $P$ by attaching $P'$ disjointly except at $v_1$ and adding edge $e'$ from $v'$ in $P'$ back to~$v_1$.
\end{enumerate}
\item Repeat step (2) for each edge $e'\neq e$ at vertex $v_1$.
\item Continue with the induction step unless $v_1=\tau(p)$. If $v_1=\tau(p)$, then output $\mathsf{Pict}(\Gamma,p)=P$.
\end{enumerate}
\end{definition}

\begin{remark}
The map $\mathsf{Pict}$ has the property that the set of all paths in $\Gamma$ from $\mathbbm{1}$ to $\tau(p)$ is in
bijection with the set of all paths in $\mathsf{Pict}(\Gamma,p)$ from $\mathbbm{1}$ to $\tau(p)$ such that
the labels of the paths are preserved.
\end{remark}

\begin{example}
\label{example.loop}
Let us compute the example $\mathsf{Pict}(\Gamma,ab\square)$, where $\Gamma = 
\mathsf{Mc}\circ \mathsf{KR}(D_2\cup \{\square\},\{a,b,\square\})$. The McCammond expansion is
given in Figure~\ref{figure.mccammond} if we attach to each vertex an edge labelled $\square$ to the ideal
$\{\square\}$. The straight path is $ab\square$. The vertex labelled $a$ in Figure~\ref{figure.mccammond}
has four dashed edges coming in. Two are labelled $a$ and two are labelled $b$. By the algorithm described
in Definition~\ref{definition.pict}, the long dashed arrows labelled $a$ and $b$ give rise to loops of length 4.
The other two dashed arrows give rise to loops of length 2. Repeating the process yields the loop graph
in Figure~\ref{figure.loop}.
\end{example}

\subsection{Kleene expressions}
\label{section.kleene}

Denote the set of all paths in a loop graph $G$ starting at $\mathbbm{1}$ and ending at $\tau(p)$, where $p$ is the 
straight line the loop graph is based on, by $\mathcal{P}_G$. We represent a path $q\in \mathcal{P}_G$ by
\[
	\mathbbm{1} \stackrel{a_1}{\longrightarrow} v_1 \stackrel{a_2}{\longrightarrow} \cdots \stackrel{a_k}
	{\longrightarrow} v_k=\tau(p),
\]
where $v_i$ are vertices in $G$ and $a_i \in A$ are the labels on the edges.

There is a simple inductive way to describe $\mathcal{P}_G$ using \defn{Kleene expressions}
(see~\cite[Section~1.3]{RhodesSchilling.2019b}).
Given a set $L$, define $L^0 = \{\varepsilon \}$ given by the empty string, $L^1 = L$, and recursively
$L^{i+1} = \{wa \mid w \in L^i, a \in L\}$ for each integer $i>0$. Then the \defn{Kleene star} is
\[
	L^\star = \bigcup_{i\geqslant 0} L^i.
\]
A Kleene expression only involves letters in $A$, concatenation, unions, and $\star$.
To obtain a Kleene expression for $\mathcal{P}_G$, perform the following doubly recursive procedure:

\smallskip
\noindent
\textbf{Algorithm 1.} Assume that the straight line path corresponding to the loop graph $G$ is indexed as
\begin{center}
\begin{tikzpicture}[auto]
\node (A) at (0, 0) {$\mathbbm{1}$};
\node (B) at (1.5,0) {$1$};
\node(C) at (3,0) {$2$};
\node(Cp) at (4.5,0) {};
\node(Dp) at (6,0) {};
\node(F) at (5.3,0) {$\cdots$};
\node(D) at (7.5,0) {$\tau(p)$};
\draw[edge,thick] (A) -- (B);
\draw[edge,thick] (B) -- (C);
\draw[edge,thick] (C) -- (Cp);
\draw[edge,thick] (Dp) -- (D);
\end{tikzpicture}
\end{center}

\smallskip
\noindent
\textbf{Induction basis:} Start at vertex $\mathbbm{1}$ and with the empty expression $E$.

\smallskip
\noindent
\textbf{Induction step:} Suppose one is at vertex $i\neq \tau(p)$ (or $\mathbbm{1}$) on the straight line path underlying $G$.
\begin{enumerate}
\item Continue to the next vertex $i+1$ (or $1$) on the straight line path underlying $G$ and append the label $a$
on the edge from $i \stackrel{a}{\longrightarrow} i+1$ (or $\mathbbm{1} \stackrel{a}{\longrightarrow} 1$) to $E$.
\item If there are loops $\ell_1,\ell_2,\ldots,\ell_k$ at vertex $i+1$ (or $1$), append the formal expression
\[
	\{\ell_1,\ell_2,\ldots,\ell_k\}^\star
\]
to $E$. The loops $\ell_1,\ell_2,\ldots,\ell_k$ are in one-to-one correspondence with the edges coming into vertex
$i+1$.
\item If $i+1\neq \tau(p)$, continue with the next induction step. Else stop and output $E$.
\end{enumerate}

\smallskip
\noindent
\textbf{Algorithm 2.} For each symbol $\ell_i$ in the expression for $E$, do the following:
\begin{enumerate}
\item Consider the loop $\ell_i = \left( v_0 \stackrel{a_1}{\longrightarrow} v_1 \stackrel{a_2}{\longrightarrow} \cdots 
\stackrel{a_k}{\longrightarrow} v_k=v_0 \right)$ from vertex $v_0$ to $v_0$ in $G$. Consider the subgraph
of $G$ with straight line $v_1 \stackrel{a_2}{\longrightarrow} \cdots \stackrel{a_k}{\longrightarrow} v_k$ and all
further loops that are attached to any of the vertices $v_i$ in $G$. Attach $\mathbbm{1}$ to $v_1$. The resulting
graph $G^{(i)}$ is a new loop graph. Perform Algorithm 1 on $G^{(i)}$ to obtain a Kleene expression $E^{(i)}$.
Replace the symbol $\ell_i$ in $E$ by $E^{(i)}$.
\item Continue this process until $E$ does not contain any further expressions $\ell_i$ for some loop $\ell_i$,  that
is, $E$ only contains unions, $\star$ and elements in the alphabet $A$.
Then the Kleene expression for $\mathcal{P}_G$ is $E$.
\end{enumerate}
The resulting expressions can be made into unionless expressions by using \defn{Zimin words}
\begin{equation}
\label{equation.zimin}
	\{a\}^\star = a^\star \qquad \text{and} \qquad \{a, b\}^\star = (a^\star b)^\star a^\star \qquad \text{for $a,b\in A$.}
\end{equation}
Expressions for larger unions can be obtained by induction using~\eqref{equation.zimin}.

\begin{example}
\label{example.kleene}
Let us continue Example~\ref{example.loop} and compute the Kleene expression for
$\mathcal{P}_{\mathsf{Pict}(\Gamma,ab\square)}$ for $\Gamma = 
\mathsf{Mc}\circ \mathsf{KR}(D_2\cup \{\square\},\{a,b,\square\})$.
By Algorithm 1, we obtain the expression 
\[
	E = a \{\ell_1,\ell_2,\ell_3,\ell_4\}^\star b \ell_5^\star \square.
\]
Using Algorithm 2 repeatedly for $\ell_1,\ldots,\ell_5$, we obtain
\[
\begin{split}
	\ell_1 &= a (b(aa)^\star b)^\star b (aa)^\star ab,\\
	\ell_2 &= a(b(aa)^\star b)^\star a,\\
	\ell_3 &= b (a(bb)^\star a)^\star a (bb)^\star ba,\\
	\ell_4 &= b(a(bb)^\star a)^\star b,\\
	\ell_5 &= a(bb)^\star a.
\end{split}
\]
\end{example}

\subsection{From Kleene expressions to rational functions}
\label{section.rational}

Our aim is to evaluate the expressions for $\Psi_w^{\mathcal{M}(S,A)}$ in Theorems~\ref{theorem.stationary}
and~\ref{theorem.stationary general}. Let $G$ be a loop graph with straight line path $p$. Define
\begin{equation}
\label{equation.psi G}
	\Psi_G(x_1,\ldots,x_n) = \sum_{q} \prod_{a \in q} x_a,
\end{equation}
where the sum is over all paths $q$ from $\mathbbm{1}$ to $\tau(p)$ in $G$.
In~\cite[Definition 1.3]{RhodesSchilling.2019b} this is also called the \defn{normal distribution} of the loop graph $G$.
Note that $\Psi_w^{\mathcal{M}(S,A)}$ is the sum of $\Psi_G(x_1,\ldots,x_n)$ for various loop graphs with straight line 
paths $p$ such that $\tau(p)=w$ (see also~\cite[Theorem 1.4]{RhodesSchilling.2019b}).
The Kleene expressions from Section~\ref{section.kleene} give us an expression
for the set of relevant paths $q$ in $G$. Now we discuss how to get from the Kleene expressions to rational functions.

The main idea is that concatenation in Kleene expressions corresponds to products, unions corresponds to sums,
and $\star$ corresponds to the geometric series. More concretely, for a path $p=a_1\cdots a_k$ we obtain
\[
	\prod_{a\in p} x_a = x_{a_1} x_{a_2} \cdots x_{a_k}.
\]
For $\star$-expressions with a single letter $a$, we obtain
\[
	\sum_{s\in a^\star} \prod_{i \in s} x_i = \sum_{\ell=0}^\infty x_a^\ell = \frac{1}{1-x_a}.
\]
Similarly
\[
	\sum_{s \in \{a,b\}^\star} \prod_{i\in s} x_i = \sum_{s \in a^\star (ba^\star)^\star} \prod_{i\in s} x_i
	= \frac{1}{1-x_a} \cdot \frac{1}{1-\frac{x_b}{1-x_a}}
	= \frac{1}{1-x_a-x_b}.
\]
In general, using the recursion~\eqref{equation.zimin} we derive by induction
\begin{equation}
\label{equation.geometric}
	\sum_{s \in \{a_1,a_2,\ldots,a_n\}^\star} \prod_{i\in s} x_i = \frac{1}{1-x_{a_1} - x_{a_2} - \cdots - x_{a_n}}.
\end{equation}

\begin{example}
\label{example.psi limit}
Let us now compute
\[
	\Psi_{ab\square} = \sum_{p\in E} \prod_{a\in p} x_a
\]
for the Kleene expression $E$ of Example~\ref{example.kleene}. We find (see also~\cite[Example 3.8]{RhodesSchilling.2019b})
\begin{equation*}
\begin{split}
	\Psi_{ab\square} &= \frac{x_a x_b x_\square}
	{\left(1-\frac{x_a^2x_b^2}{\left(1-\frac{x_b^2}{1-x_a^2}\right)(1-x_a^2)} -\frac{x_a^2}{1-\frac{x_b^2}{1-x_a^2}}
	-\frac{x_a^2x_b^2}{\left(1-\frac{x_a^2}{1-x_b^2}\right)(1-x_b^2)} -\frac{x_b^2}{1-\frac{x_a^2}{1-x_b^2}}\right)
	\left(1-\frac{x_a^2}{1-x_b^2}\right)}\\
	&=\frac{x_a x_b x_\square(1-x_b^2)}
         {\left(1-\frac{2x_a^2x_b^2}{1-x_a^2-x_b^2} -\frac{x_a^2(1-x_a^2)}{1-x_a^2-x_b^2}
	-\frac{x_b^2(1-x_b^2)}{1-x_a^2-x_b^2}\right) (1-x_a^2-x_b^2)}\\
	&=\frac{x_a x_b x_\square(1-x_b^2)}
         {1-2x_a^2-2x_b^2+(x_a^2-x_b^2)^2}.
\end{split}
\end{equation*}
In the limit as $x_\square \to 0$, we obtain
\[
	\lim_{x_\square \to 0} \Psi_{ab\square} = \frac{1-x_b^2}{8}.
\]
In~\cite[Example 3.8]{RhodesSchilling.2019b}, the remaining stationary distributions were computed 
using that $x_a+x_b+x_\square=1$ and by taking the limit $x_\square \to 0$
\begin{equation*}
\begin{split}
	\Psi_\square & = x_\square \qquad \qquad \qquad \qquad \qquad \qquad \;\;
	 \stackrel{x_\square\to 0}{\longrightarrow} \qquad 0\\
	\Psi_{a\square} &= \frac{x_a(1-x_a^2-x_b^2)x_\square}{1-2x_a^2-2x_b^2+(x_a^2-x_b^2)^2}
	\qquad \stackrel{x_\square\to 0}{\longrightarrow} \qquad \frac{x_a}{4}\\
	\Psi_{ab\square} &= \frac{x_a x_b x_\square(1-x_b^2)} {1-2x_a^2-2x_b^2+(x_a^2-x_b^2)^2} 
	\qquad \stackrel{x_\square\to 0}{\longrightarrow} \qquad \frac{1-x_b^2}{8}\\
	\Psi_{aba\square} &= \frac{x_a^2 x_b x_\square}{1-2x_a^2-2x_b^2+(x_a^2-x_b^2)^2}
	\qquad \stackrel{x_\square\to 0}{\longrightarrow} \qquad \frac{x_a}{8}\\
	\Psi_{abab\square} &= \frac{x_a^2 x^2_b x_\square}{1-2x_a^2-2x_b^2+(x_a^2-x_b^2)^2}
	\qquad \stackrel{x_\square\to 0}{\longrightarrow} \qquad \frac{x_a x_b}{8}\\
	\Psi_{a^2\square} &= \frac{x_a^2 (1-x_a^2) x_\square}{1-2x_a^2-2x_b^2+(x_a^2-x_b^2)^2}
	\qquad \stackrel{x_\square\to 0}{\longrightarrow} \qquad \frac{x_a (1+x_a)}{8}\\
	\Psi_{a^2b\square} &= \frac{x_a^2 x_b x_\square}{1-2x_a^2-2x_b^2+(x_a^2-x_b^2)^2}
	\qquad \stackrel{x_\square\to 0}{\longrightarrow} \qquad \frac{x_a}{8}\\
	\Psi_{a^2ba\square} &= \frac{x_a^3 x_b x_\square}{1-2x_a^2-2x_b^2+(x_a^2-x_b^2)^2}
	\qquad \stackrel{x_\square\to 0}{\longrightarrow} \qquad \frac{x_a^2}{8}
\end{split}
\end{equation*}
and similarly for the cases with $a$ and $b$ interchanged by symmetry.

For a word $w$ in $\{a,b\}$, denote by $[w]$ the corresponding element in $D_2$. For example, in $D_2$ we have 
$[a]=[bab]=[b^2a]$. Note that by Theorem~\ref{theorem.stationary general}
\[
	\Psi^{\mathcal{M}(D_2,\{a,b\})}_s=\frac{1}{4} \qquad \text{for all $s\in D_2$}
\]
by summing the appropriate results for $\lim_{x_\square \to 0} \Psi_{w\square}$ as above. For example,
\[
	 \Psi^{\mathcal{M}(D_2,\{a,b\})}_{[ab]} = \lim_{x_\square \to 0} \left(\Psi_{ab\square} + \Psi_{ba\square}
	 +\Psi_{a^2ba\square} + \Psi_{b^2ab\square} \right)
	 = \frac{1-x_b^2}{8} + \frac{1-x_a^2}{8} + \frac{x_a^2}{8} + \frac{x_b^2}{8} = \frac{1}{4}.
\]
This shows that the stationary distribution is uniform.
\end{example}

\begin{remark}
Note that the Markov chain in~\eqref{equation.transition diagram} is not ergodic since the greatest common
divisor of the cycle length is 2 and not 1. We can make the Markov chain ergodic by introducing a new generator $c$,
which acts as the identity. In other words, this would introduce loops at each vertex in~\eqref{equation.transition diagram}
labeled $c$. In turn, this would introduce loops labeled $c$ at each vertex in the McCammond expansion
in Figure~\ref{figure.mccammond} and the loop graph in Figure~\ref{figure.loop}. This would change the Kleene
expressions in Example~\ref{example.kleene} to
\[
	E = a \{\ell_1,\ell_2,\ell_3,\ell_4,c\}^\star b \{\ell_5,c\}^\star \square
\]
with
\[
\begin{split}
	\ell_1 &= a \{b\{ac^\star a,c\}^\star b,c\}^\star b \{ac^\star a,c\}^\star a c^\star b,\\
	\ell_2 &= a\{b\{ac^\star a,c\}^\star b,c\}^\star a,\\
	\ell_3 &= b \{a\{b c^\star b,c\}^\star a,c\}^\star a \{bc^\star b,c\}^\star b c^\star a,\\
	\ell_4 &= b\{a\{bc^\star b,c\}^\star a,c\}^\star b,\\
	\ell_5 &= a\{bc^\star b,c\}^\star a.
\end{split}
\]
In this setting, we find
\begin{multline*}
	\Psi_{ab\square} =
	\frac{x_a x_b x_\square}
	{\left(1
	-\frac{x_a^2x_b^2}{\left(1-\frac{x_b^2}{1-\frac{x_a^2}{1-x_c}-x_c}-x_c\right)\left(1-\frac{x_a^2}{1-x_c}-x_c\right)(1-x_c)} 
	-\frac{x_a^2}{1-\frac{x_b^2}{1-\frac{x_a^2}{1-x_c}-x_c}-x_c}\right.}\\
	\frac{1}{\left.
	-\frac{x_a^2x_b^2}{\left(1-\frac{x_a^2}{1-\frac{x_b^2}{1-x_c}-x_c}-x_c\right)\left(1-\frac{x_b^2}{1-x_c}-x_c\right)(1-x_c)} 
	-\frac{x_b^2}{1-\frac{x_a^2}{1-\frac{x_b^2}{1-x_c}-x_c}-x_c}-x_c\right) 
	 \left(1-\frac{x_a^2}{1-\frac{x_b^2}{1-x_c}-x_c}-x_c\right)}\\
	 = \frac{x_a x_b x_\square (1+x_b-x_c)(1-x_b-x_c)}{(1-x_a-x_b-x_c)(1+x_a+x_b-x_c)(1-x_a+x_b-x_c)(1+x_a-x_b-x_c)}.
\end{multline*}
Using $x_a+x_b+x_c+x_\square=1$, we see that the term $(1-x_a-x_b-x_c)$ in the denominator cancels with
$x_\square$ in the numerator. Hence in the limit $x_\square \to 0$, we obtain
\[
	\lim_{x_\square \to 0} \Psi_{ab\square} = \frac{(1+x_b-x_c)(1-x_b-x_c)}{8(x_a+x_b)}.
\]
As $x_c\to 0$, we recover the result from Example~\ref{example.psi limit}.
\end{remark}

\subsection{Mixing time}
\label{section.mixing}

Recall from Section~\ref{section.markov}, that for a given small $\epsilon>0$, the \defn{mixing time} $t_\mathsf{mix}$ 
is the smallest $t$ such that
\[
	\| T^t \nu - \Psi \| \leqslant \epsilon.
\]
Many references about mixing time can be found in~\cite{Diaconis.2011,LevinPeres.2017}.

Let $\tau$ be the first time that the Markov chain hits the ideal (when starting at $\mathbbm{1}$ in
$\mathsf{RCay}(S,A)$ or $\mathsf{Mc}\circ \mathsf{KR}(S,A)$). Denote by $\mathsf{Pr}(\tau> t)$ the probability
that $\tau$ is bigger than a given $t$. In~\cite{ASST.2015}, it was shown that $\mathsf{Pr}(\tau>t)$ gives a bound
on the mixing time.

\begin{theorem} \cite{ASST.2015}
\label{theorem.ASST}
Let $S$ be a finite semigroup whose minimal ideal $K(S)$ is a left zero semigroup and let $T$ be the transition
matrix of the associated Markov chain. Then
\[
	\| T^t \nu - \Psi \| \leqslant \mathsf{Pr}(\tau>t).
\]
\end{theorem}

In~\cite{RhodesSchilling.2020}, we provide a way to compute $\mathsf{Pr}(\tau>t)$ from particular rational
expressions for the stationary distribution. Let $G$ be a loop graph and recall $\Psi_G(x_1,\ldots,x_n)$ 
from~\eqref{equation.psi G}, which is a rational function in $x_1,\ldots,x_n$. Let $\mathsf{Pr}_G(\tau\geqslant t)$ be the 
probability that the length of the paths in the loop graph $G$ from $\mathbbm{1}$ to $s$ in the ideal is weakly bigger than $t$.
Let $\Psi^{\geqslant t}_G(x_1,\ldots,x_n)$ be the truncation of the formal power series associated to the rational function
$\Psi_G(x_1,\ldots,x_n)$ to terms of degree weakly bigger than $t$ and let $\Psi_G^{<t}(x_1,\ldots,x_n)$ be the truncation
of the formal power series associated to the rational function $\Psi_G(x_1,\ldots,x_n)$ to terms of degree strictly smaller
than $t$. Note that
\[
	\Psi_G(x_1,\ldots,x_n) = \Psi_G^{<t}(x_1,\ldots,x_n) + \Psi^{\geqslant t}_G(x_1,\ldots,x_n).
\]

\begin{theorem} \cite{RhodesSchilling.2020}
\label{theorem.main}
Suppose the Markov chain satisfies the conditions of Theorem~\ref{theorem.ASST}. 
If $\Psi_G(x_1,\ldots,x_n)$ is represented by a rational function such that each term of degree $\ell$ 
in its formal power sum expansion corresponds to a path in $G$ of length $\ell$, we have
\[
	\mathsf{Pr}_G(\tau \geqslant t) = \frac{\Psi^{\geqslant t}_G(x_1,\ldots,x_n)}{\Psi_G(x_1,\ldots,x_n)}
	= 1 - \frac{\Psi_G^{<t}(x_1,\ldots,x_n)}{\Psi_G(x_1,\ldots,x_n)}.
\]
\end{theorem}

By Markov's inequality (see for example~\cite{LevinPeres.2017,DevroyeLugosi.2001}), we have
\begin{equation}
\label{equation.Markov inequality}
	\mathsf{Pr}(\tau>t) \leqslant \frac{E[\tau]}{t+1},
\end{equation}
where $E[\tau]$ is the expected value for $\tau$, the first time the walk hits the ideal. Hence knowing $E[\tau]$
gives an upper bound on the mixing time. In~\cite{RhodesSchilling.2020}, we find a way to compute $E[\tau]$
from certain representations of the stationary distribution.

\begin{theorem} \cite{RhodesSchilling.2020}
\label{theorem.main1}
Suppose the Markov chain satisfies the conditions of Theorem~\ref{theorem.ASST}.
If $\Psi_G(x_1,\ldots,x_n)$ is represented by a rational function such that each term of degree $\ell$ in its formal power 
sum expansion corresponds to a path in $G$ of length $\ell$, we have
\begin{equation}
\label{equation.EG}
	E_G[\tau] = \left( \sum_{i=1}^n x_i \frac{\partial}{\partial x_i} \right) \ln \Psi_G(x_1,\ldots,x_n).
\end{equation}
\end{theorem}

\subsection{Another example}
\label{section.example}

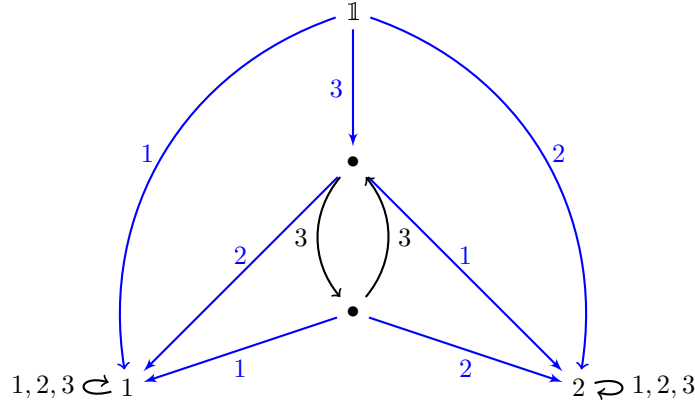
\begin{figure}[t]
\begin{center}
\begin{tikzpicture}[auto]
\node (A) at (0, 0) {$\mathbbm{1}$};
\node (B) at (-3,-5) {$1$};
\node(C) at (3,-5) {$2$};
\node(D) at (0,-2) {$\bullet$};
\node(E) at (0,-4) {$\bullet$};
\path (A) edge[->,thick, blue, bend right = 40] node[midway,left] {$1$} (B)
              edge[->,thick, blue, bend left = 40] node[midway,right] {$2$} (C);
\draw[edge,thick,blue] (A) -- (D) node[midway, left] {$3$};
\draw[edge,thick,blue] (D) -- (B) node[midway, above] {$2$};
\draw[edge,thick,blue] (D) -- (C) node[midway, above] {$1$};
\draw[edge,thick,blue] (E) -- (B) node[midway, below] {$1$};
\draw[edge,thick,blue] (E) -- (C) node[midway, below] {$2$};
\path (D) edge[->,thick, bend right = 40] node[midway,left] {$3$} (E);
\path (E) edge[->,thick, bend right = 40] node[midway,right] {$3$} (D);
\path (B) edge[->,thick, loop left] node {$1,2,3$} (B);
\path (C) edge[->,thick, loop right] node {$1,2,3$} (C);
\end{tikzpicture}
\end{center}
\caption{\label{figure.right cayley linear} The right Cayley graph $\mathsf{RCay}(S,A)$ of the semigroup that
gives the Markov chain in Section~\ref{section.example}.}
\end{figure}

Let us illustrate the concepts and algorithms in terms of another example. Consider the Markov chain with state
space $\Omega=\{\mathbf{1},\mathbf{2}\}$ given by the transition diagram:
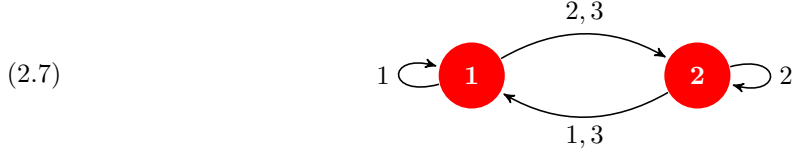
\begin{equation}
\label{equation.markov linear}
\raisebox{-1cm}{
\begin{tikzpicture}[->,>=stealth',shorten >=1pt,auto,node distance=3cm,
                    semithick]
  \tikzstyle{every state}=[fill=red,draw=none,text=white]

  \node[state]         (I) {$\mathbf{1}$};
  \node[state]         (A) [right of=I] {$\mathbf{2}$};
  \path (I) edge [bend left] node {$2,3$} (A)
               edge [loop left] node {$1$} (I)
           (A) edge[bend left] node {$1,3$} (I)
                 edge [loop right] node {$2$} (A);
\end{tikzpicture}}
\end{equation}
The transition matrix of this Markov chain is given by
\[
	T = \begin{pmatrix}
	x_1 & x_1+x_3\\
	x_2+x_3 & x_2
	\end{pmatrix}.
\]
Pick as generators of the semigroup $A=\{1,2,3\}$. Then the right Cayley graph of the semigroup that gives
the above Markov chain by left multiplication is depicted in Figure~\ref{figure.right cayley linear}, where $K(S)=\{1,2\}$.
Indeed, the left action $11=1$, $21=2$ and $31=2$, which gives all the edges out of $\mathbf{1}$ 
in~\eqref{equation.markov linear}.
Similarly, the left action $12=1$, $22=2$, and $32=1$, which gives all the edges out of $\mathbf{2}$ 
in~\eqref{equation.markov linear}.

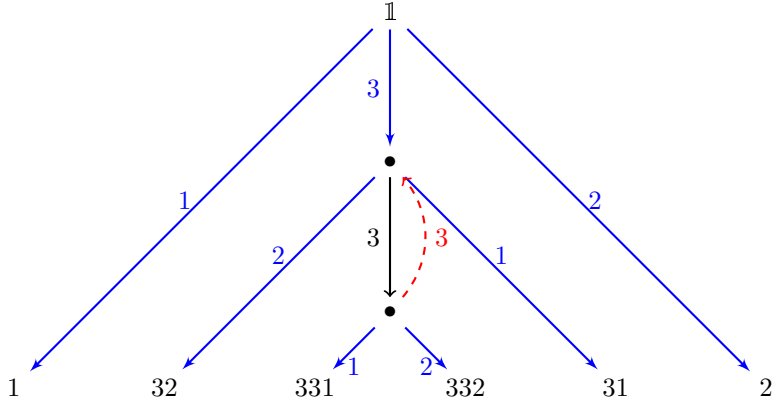
\begin{figure}[t]
\begin{center}
\begin{tikzpicture}[auto]
\node (A) at (0, 0) {$\mathbbm{1}$};
\node (B) at (-5,-5) {$1$};
\node(C) at (5,-5) {$2$};
\node(D) at (0,-2) {$\bullet$};
\node(E) at (0,-4) {$\bullet$};
\node(F) at (-3,-5) {$32$};
\node(G) at (-1,-5) {$331$};
\node(H) at (1,-5) {$332$};
\node(I) at (3,-5) {$31$};
\draw[edge,thick,blue] (A) -- (B) node[midway, left] {$1$};
\draw[edge,thick,blue] (A) -- (C) node[midway, right] {$2$};
\draw[edge,thick,blue] (A) -- (D) node[midway, left] {$3$};
\draw[edge,thick,blue] (D) -- (F) node[midway, above] {$2$};
\draw[edge,thick,blue] (D) -- (I) node[midway, above] {$1$};
\draw[edge,thick,blue] (E) -- (G) node[midway, below] {$1$};
\draw[edge,thick,blue] (E) -- (H) node[midway, below] {$2$};
\path (D) edge[->,thick] node[midway,left] {$3$} (E);
\path (E) edge[->,thick, red, dashed, bend right = 40] node[midway,right] {$3$} (D);
\end{tikzpicture}
\end{center}
\caption{\label{figure.mccammond linear} $\mathsf{Mc}\circ \mathsf{KR}(S,A)$ of $\mathsf{RCay}(S,A)$ in
Figure~\ref{figure.right cayley linear} with loops on the ideal omitted.}
\end{figure}

The McCammond and Karnofsky--Rhodes expansion of the right Cayley graph is given in 
Figure~\ref{figure.mccammond linear}.
The Kleene expression for all paths from $\mathbbm{1}$ to $32 \in K(S)$ is $3(33)^\star 2$ and similarly for
paths with other endpoints. From this we easily compute
\[
\begin{aligned}
	\Psi_1 &= x_1 & \qquad \Psi_2&=x_2,\\
	\Psi_{32} &= \frac{x_2 x_3}{1-x_3^2} & \qquad \Psi_{31} &= \frac{x_1 x_3}{1-x_3^2},\\
	\Psi_{331} &= \frac{x_1 x_3^2}{1-x_3^2} & \qquad \Psi_{332} &= \frac{x_2 x_3^2}{1-x_3^2}.
\end{aligned}
\]
Furthermore,
\[
\begin{split}
	\Psi^{\mathcal{M}(S,A)}_1 &= \Psi_1 + \Psi_{32} + \Psi_{331} = \frac{x_1+x_2 x_3}{1-x_3^2},\\
	\Psi^{\mathcal{M}(S,A)}_2 &= \Psi_2 + \Psi_{31} + \Psi_{332} = \frac{x_2+x_1 x_3}{1-x_3^2}.
\end{split}
\]
Using that $x_1+x_2+x_3=1$, we find indeed that $\Psi^{\mathcal{M}(S,A)}_1+\Psi^{\mathcal{M}(S,A)}_2=1$.

Since the expressions for $\Psi^{\mathcal{M}(S,A)}_1$ and $\Psi^{\mathcal{M}(S,A)}_2$ were computed
directly from $\mathsf{Mc} \circ \mathsf{KR}(S,A)$ (or the corresponding loop graphs) without using that 
$x_1+x_2+x_3=1$, each term of degree $\ell$ in the expansion of the rational function corresponds to a path
in the graph. Hence we may use Theorem~\ref{theorem.main1} to give an upper bound on the mixing time
\[
	E_1[\tau] = \left( x_1 \frac{\partial}{\partial x_1} +  x_2 \frac{\partial}{\partial x_2} + x_3 \frac{\partial}{\partial x_3} \right)
	\ln \Psi^{\mathcal{M}(S,A)}_1
	= \frac{x_1}{x_1+x_2 x_3} + \frac{2 x_2 x_3}{x_1+x_2 x_3} + \frac{2 x_3^2}{1-x_3^2}.
\]
Inserting $x_1=x_2=x_3=\frac{1}{3}$ yields $E_1[\tau]=E_2[\tau]=\frac{3}{2}$, so that $t_{\mathsf{mix}} \leqslant 3$ if
$\epsilon = \frac{1}{2}$.

\bibliographystyle{alpha}

\end{document}